\documentstyle[amsfonts,12pt]{article}
\newtheorem{th}{Theorem}[section]
\newtheorem{lem}[th]{Lemma}

\newtheorem{cor}[th]{Corollary}
\newtheorem{defn}[th]{Definition}
\newenvironment{defn-new}{\begin{defn} \em}{\end{defn}}
\newtheorem{rem}[th]{Remark}
\newenvironment{rem-new}{\begin{rem} \em}{\end{rem}}
\newtheorem{ex}[th]{Example}
\newenvironment{ex-new}{\begin{ex} \em}{\end{ex}}

\newenvironment{notation-new}{\begin{rem} \em}{\end{rem}}

\newenvironment{agr-new}{\begin{rem} \em}{\end{rem}}

\makeatletter \@addtoreset{equation}{section} \makeatother

\makeatletter \@addtoreset{figure}{section} \makeatother

\begin{document}

\begin{center}
{\Large Semi-invariant Conformal submersions with horizontal Reeb vector field }

\bigskip {\large Uday Chand De, Shashikant Pandey and Punam Gupta }
\end{center}


\noindent{\bf Abstract.}
The present paper deals with the characterization of a new submersion named
semi-invariant conformal $\zeta ^{\perp }$-Riemannian submersion from almost
contact metric manifolds onto Riemannian manifolds which is the
generalization of some known submersions on Riemannian manifolds. We give
important and adequate conditions for such submersions to be totally
geodesic and harmonic. Also, few examples are examined for such submersions
endowed with horizontal Reeb vector field.

%
%
%
%
%
%
%
%
%



\section{Introduction}

In Riemannian geometry, to develop Riemannian manifolds with non-negative
sectional curvature is an iconic problem. For this, O'Neill \cite{Neill} and
Gray \cite{Gray} initiate to study Riemannian submersions connecting
Riemannian manifolds. Riemannian submersions have several applications in
mathematical physics, specially in the Yang--Mills theory \cite{Bourg,
Watson}, Kaluza-Klein theory \cite{Bourg-1,Ianus}, supergravity and
superstring theories \cite{Ianus-Visi, Mustafa}. Submersions were studied by
many authors between differential manifolds like almost Hermitian manifolds,
K\"{a}hler manifold, almost contact manifolds, nearly $K$-cosymplectic
manifold, Sasakian manifold etc. The authors studied semi-Riemannian
submersion and Lorentzian submersion \cite{Falci}, anti-invariant submersion
\cite{Kupeli-Murathan, Murathan, Sahin, Tan}, anti-invariant $\zeta ^{\perp }
$-submersion \cite{Lee}, semi-invariant submersion \cite{Sahin-1}, \cite%
{RKRS}, semi-invariant $\zeta ^{\perp }$-submersion \cite{Akyol-Sari-Eksoy}
and many other type of submersions studied by several authors \cite{Chen,
Kupeli,Sahin-2, Watson-1,Ianus-1,Park, Park-T}, hemi-slant submersion \cite%
{Tastan}, anti-holomorphic semi-invariant submersion \cite{Tastan-14}\ etc.
between differential manifolds with different structures.

In comparison of conformal submersions, Riemannian submersions are specific.
The conformal maps don't preserve distance between points but they preserve
angle between vector fields which allows us to transfer certain properties
of manifolds to another manifolds by deforming such properties. Ornea \cite%
{Ornea} initiated the theory of conformal submersions between Riemannian
manifolds. Later, it was studied by many authors \cite%
{Akyol,Akyol-Sahin-16,Akyol-Sahin-17, Akyol-Sahin,Gundu-Akyol}.

In the present manuscript, we establish the definition of semi-invariant
conformal $\zeta ^{\perp}$-Riemannian submersion (in short, sic $\zeta
^{\perp }$-Rs) from almost contact metric (in short, acm) manifolds onto
Riemannian manifolds, which is the generalization of conformal
anti-invariant $\zeta ^{\perp }$-submersions, semi-invariant conformal
submersions and many others. In section \ref{sect-pre}, we give
preliminaries about Sasakian manifolds, Riemannian submersions and conformal
submersions. In section \ref{sect-con}, we define sic $\zeta ^{\perp }$-Rs
from acm manifolds onto Riemannian manifolds with examples. We give the
necessary and sufficient (briefly, ns) conditions for distributions to be
integrable. After that, we establish the ns condition for sic $\zeta ^{\perp
}$-Rs to be homothetic map. Then we investigate the geometry of leaves of
horizontal and vertical distributions and obtain ns condition for the
distribution to be absolutely geodesic foliation on Sasakian manifolds. At
last, we obtain a condition when the fibers of sic $\zeta ^{\perp }$-Rs from
Sasakian manifolds onto Riemannian manifolds to be locally product manifold.
In section \ref{sect-har}, we provide ns conditions for a sic $\zeta ^{\perp
}$-Rs to be harmonic. We also investigate the ns conditions for such
submersions to be completely geodesic.

\section{Preliminaries\label{sect-pre}}

In this segment, we recall some definitions which will be required to study
throughout the paper. \newline
Let $(M^{\ast},\varphi_{\ast} ,\zeta ,\eta_{\ast} ,g_{m})$ be an odd
dimensional acm manifold \cite{Blair} consisting of a $(1,1)$-tensor field $%
\varphi_{\ast} $, a vector field $\zeta $, a $1$-form $\eta_{\ast} $ and a
Riemannian metric $g_{m}$. Then
\begin{equation}
\varphi_{\ast} ^{2}P=-P+\eta_{\ast} (P)\zeta ,  \label{EQ2.1}
\end{equation}%
\begin{equation}
\varphi_{\ast} \,\zeta =0,\ \eta_{\ast} \circ \varphi_{\ast} =0,\
\eta_{\ast} (\zeta )=1,\ g_{m}(P,\zeta )=\eta_{\ast} (P),  \label{EQ2.2}
\end{equation}%
and
\begin{equation}
g_{m}(\varphi_{\ast} P,\varphi_{\ast} Q)=g_{m}(P,Q)-\eta_{\ast}
(P)\eta_{\ast} (Q),\quad g_{m}(\varphi_{\ast} P,Q)=-g_{m}(P,\varphi_{\ast}
Q),  \label{EQ2.3}
\end{equation}%
for any vector fields $P,Q\in \Gamma (TM^{\ast})$.\newline

An acm manifold $(M^{\ast},\varphi_{\ast} ,\zeta ,\eta_{\ast} ,g_{m})$ is
called \emph{Sasakian }\cite{Sasaki}\emph{\ }if
\begin{equation}
(\nabla _{P}\varphi_{\ast} )Q=g_{m}(P,Q)\zeta -\eta_{\ast} (Q)P
\label{EQ2.4}
\end{equation}%
and
\begin{equation}
\nabla _{P}\zeta =-\varphi_{\ast} P  \label{EQ2.5}
\end{equation}%
for any vector fields $P,Q\in \Gamma (TM^{\ast})$, where $\nabla $ is the
Levi-Civita connection.

Now, we give some useful definitions of submersions to extend our definition
of semi-invariant conformal Riemannian submersions endowed with the
horizontal Reeb vector field.

\begin{defn}
\textrm{{\cite{Neill}}}\ Let $(M^{\ast},g_{m})$ and $(N^{\ast},g_{n})$ be
Riemannian manifolds, where $\dim (M^{\ast})=m$, $\dim (N^{\ast})=n$ and $m>n
$. A Riemannian submersion $\Psi :M^{\ast}\rightarrow N^{\ast}$ is a map of $%
M^{\ast}$ onto $N^{\ast}$ satisfying the following conditions:

\begin{description}
\item[(i)] $\Psi$ has maximal rank.

\item[(ii)] The metric of horizontal vectors is preserved by the
differential $\Psi _{\ast }$.
\end{description}
\end{defn}

For each $z\in N^{\ast}$, $\Psi ^{-1}(z)$ is an $(m-n)$-dimensional
submanifold of $M^{\ast} $. The submanifolds $\Psi ^{-1}(z)$, $z\in N^{\ast}$%
, are called \emph{fibers}. A vector field on $M^{\ast}$ is called vertical
and horizontal if it is always tangent to fibers and orthogonal to fibers,
respectively. A vector field $P$ on $M$ is called basic if $P$ is horizontal
and $\Psi $-related to a vector field $P^{\prime }$ on $N$, that is, $\Psi
_{\ast }P_{x}=P_{\Psi _{\ast }(x)}^{\prime }$ for all $x\in M.$ The
projection morphisms on the distributions $\ker \Psi _{\ast }$ and $(\ker
\Psi _{\ast })^{\perp }$ are denoted by ${\mathcal{V}_{1}}$ and ${\mathcal{H}%
_{1}}$, respectively. The sections of ${\mathcal{V}_{1}} $ and ${\mathcal{H}%
_{1}}$ are called the vertical vector fields and horizontal vector fields,
respectively. So
\[
{\mathcal{V}_{1}}_{p}=T_{p}\left( \Psi ^{-1}(q)\right) ,\qquad {\mathcal{H}%
_{1}}_{p}=T_{p}\left( \Psi ^{-1}(q)\right) ^{\perp }.
\]

\begin{defn}
\textrm{{\cite{Ornea}}}\ Let $(M^{\ast},g_{m})$ and $(N^{\ast},g_{n})$ be
Riemannian manifolds, where $\dim (M^{\ast})=m$, $\dim (N^{\ast})=n$ and $m>n
$. A (horizontally) conformal submersion $\Psi :M^{\ast}\rightarrow N^{\ast}$
is a map of $M^{\ast}$ onto $N^{\ast}$ satisfying the following axioms:

\begin{description}
\item[(i)] $\Psi$ has maximal rank.

\item[(ii)] The angle between the horizontal vectors is preserved by the
differential $\Psi _{\ast }$,
\[
g_{n}(\Psi _{\ast }U,\Psi _{\ast }V)=\Lambda_{\ast} (p)g_{m}(U,V),\quad
U,V\in {\mathcal{H}_{1}}_{p_{\ast}}
\]%
where $\Lambda_{\ast} (p)$ is a non-zero number and $p\in M^{\ast}$. The
number $\Lambda_{\ast} (p)$ is called the square dilation of $\Psi $ at $p$,
it is necessarily non-negative. Its square root $\lambda_{\ast} (p)=\sqrt{%
\Lambda_{\ast} (p)}$ is called the dilation of $\Psi $ at $p$.
\end{description}
\end{defn}

Clearly, Riemannian submersion is a horizontally conformal submersion with $%
\lambda_{\ast} =1$.

\begin{defn}
\cite{Baird} Let $\Psi :(M^{\ast },g_{m})\rightarrow (N^{\ast },g_{n})$ be a
conformal submersion. A vector field $E_{1}$ on $M^{\ast }$ is called
projectiable if there exist a vector field $\widehat{E_{1}}$ on $N^{\ast }$
such that $\Psi _{\ast }({E_{1}}_{p})=\widehat{E_{1}}_{\Psi _{\ast }(p)}$
for any $p\in M^{\ast }.$ In this case, $E_{1}$ and $\widehat{E_{1}}$ are
called $\Psi $-related. A horizontal vector field $V$ on $M^{\ast }$ is
called basic, if it is projectiable.
\end{defn}

It is to be noted that if $\widehat{W}$ is a vector field on $N^{\ast },$
then there exists a unique basic vector field $W$ on $M^{\ast }$ which is
called the horizontal lift of $\widehat{W}.$

\begin{defn}
\cite{Baird} A horizontally conformal submersion $\Psi :(M^{\ast
},g_{m})\rightarrow (N^{\ast },g_{n})$ is called horizontally homothetic if
the slope of its dilation $\lambda _{\ast }$ is vertical, that is,
\[
{\mathcal{H}_{1}}(grad\lambda _{\ast })=0
\]%
at $p\in M^{\ast }$, where $\mathcal{H}{_{1}}$ is the projection on the
horizontal space $\left( \ker \Psi _{\ast }\right) ^{\perp }$.
\end{defn}

The second fundamental tensors of all fibers $\Psi ^{-1}(q),\ q\in N^{\ast}$
give rise to tensor field $T$ and $A$ on $M^{\ast}$ defined by O'Neill \cite%
{Neill} for arbitrary vector field $E$ and $F$, which is
\begin{equation}
T_{E}F={\mathcal{H}_{1}}\nabla _{{\mathcal{V}_{1}}E}^{M}{\mathcal{V}_{1}}F+{%
\mathcal{V}_{1}}\nabla _{{\mathcal{V}_{1}}E}^{M}{\mathcal{H}_{1}}F,
\label{EQ2.9}
\end{equation}%
\begin{equation}
A_{E}F={\mathcal{H}_{1}}\nabla _{{\mathcal{H}_{1}}E}^{M}{\mathcal{V}_{1}}F+{%
\mathcal{V}_{1}}\nabla _{{\mathcal{H}_{1}}E}^{M}{\mathcal{H}_{1}}F,
\label{EQ2.8}
\end{equation}%
where ${\mathcal{V}_{1}}$ and ${\mathcal{H}_{1}}$ are the vertical and
horizontal projections, respectively.

From equations (\ref{EQ2.9}) and (\ref{EQ2.8}), we obtain the following
equations
\begin{equation}
\nabla _{V}W=T_{V}W+\widehat{\nabla }_{V}W,  \label{EQ2.10}
\end{equation}%
\begin{equation}
\nabla _{V}P={\mathcal{H}_{1}}\nabla _{V}P+T_{V}P,  \label{EQ2.11}
\end{equation}%
\begin{equation}
\nabla _{P}V=A_{P}V+{\mathcal{V}_{1}}\nabla _{P}V,  \label{EQ2.12}
\end{equation}%
\begin{equation}
\nabla _{P}Q={\mathcal{H}_{1}}\nabla _{P}Q+A_{P}Q,  \label{EQ2.13}
\end{equation}%
for all $V,W\in \Gamma (\ker \Psi _{\ast })$ and $P,Q\in \Gamma (\ker \Psi
_{\ast })^{\perp },$ where ${\mathcal{V}_{1}}\nabla _{V}W=\widehat{\nabla }%
_{V}W.$ If $P$ is the basic vector field, then $A_{P}Q={\mathcal{H}_{1}}%
\nabla _{Q}P.$

Clearly, for $p\in M,$ $U\in {\mathcal{V}_{1}}_{p}$ and $P\in {\mathcal{H}%
_{1}}_{p}$, the linear operators
\[
{T}_{U},{A}_{P}:T_{p}M^{\ast }\rightarrow T_{p}M^{\ast }
\]%
are skew-symmetric, that is,
\begin{equation}
g_{m}({A}_{P}E,F)=-g_{m}(E,{A}_{P}F){\rm and} g_{m}({T}_{U}E,F)=-g_{m}(E,{T}_{U}F),  \label{EQ2.14}
\end{equation}%
for all $E,F\in $ $T_{p}M^{\ast }.$ It is clear, the restriction of ${T}$ to
the vertical distribution ${T}|_{\ker \Psi _{\ast }\times \ker \Psi _{\ast }}
$ is certainly the second fundamental form of the fibres of $\Psi $. From (%
\ref{EQ2.14}), it is easy to say that $\Psi $ has totally geodesic fibres if
and only if ${T}\equiv 0$.

Consider the smooth map $\Psi :(M^{\ast},g_{m})\rightarrow (N^{\ast},g_{n})$
between Riemannian manifolds. Then the differential $\Psi _{\ast }$ of $\Psi
$ can be noticed as a section of the bundle $Hom(TM^{\ast},\Psi
^{-1}TN^{\ast})\rightarrow M^{\ast}$, where $\Psi ^{-1}TN^{\ast}$ is the
bundle which has fibres $\left( \Psi ^{-1}TN^{\ast}\right)
_{x}=T_{f(x)}N^{\ast}$ , $x\in M^{\ast}$. $Hom(TM^{\ast},\Psi ^{-1}TN^{\ast})
$ has a connection $\nabla $ induced from the Riemannian connection $\nabla
^{M^{\ast}}$ and the pullback connection. The second fundamental form of $%
\Psi $ is
\begin{equation}
(\nabla \Psi _{\ast })(E,F)=\nabla _{E}^{N^{\ast}}\Psi _{\ast }F-\Psi _{\ast
}(\nabla _{E}^{M^{\ast}}F),{ \rm for all} E,F\in \Gamma (TM^{\ast}),
\label{EQ2.15}
\end{equation}%
where $\nabla ^{N^{\ast}}$ is the pullback connection (\cite{Baird,Bejancu}%
). The map $\Psi $ is said to be totally geodesic \cite{Baird} if $(\nabla
\Psi _{\ast })(E,F)=0,$ for all $E,F\in \Gamma (TM^{\ast})$.

\begin{lem}
\label{LE2}  \cite[Lemma 1.16, pp.129]{Urakawa} Let $\Psi
:(M^{\ast},g_{m})\rightarrow (N^{\ast},g_{n})$ be a smooth map between
Riemannian manifolds $(M^{\ast},g_{m})$ and $(N^{\ast},g_{n})$. Then
\begin{equation}
\Psi _{\ast }([P,Q])=\nabla _{P}^{N^{\ast}}\Psi _{\ast }Q-\nabla
_{Q}^{N^{\ast}}\Psi _{\ast }P,  \label{EQ2.16}
\end{equation}
for all $P,Q\in \Gamma (TM^{\ast})$.
\end{lem}

From Lemma \ref{LE2}, we come to the conclusion
\begin{equation}
\lbrack P,Q]\in \Gamma (ker\Psi _{\ast })  \label{EQ2.17}
\end{equation}%
for $P\in \Gamma (ker\Psi _{\ast })^{\perp }$ and $Q\in \Gamma (ker\Psi
_{\ast })$.

A smooth map $\Psi :(M^{\ast },g_{m})\rightarrow (N^{\ast },g_{n})$ is said
to be harmonic \cite{Baird} if and only if $\mathrm{trace}(\nabla \Psi
_{\ast })=0$. The tension field of $\Psi $ is the section $\tau (\Psi )$ of $%
\Gamma (\Psi ^{-1}TN)$ and defined by
\begin{equation}
\tau (\Psi )={div}\Psi _{\ast }=\sum_{i=1}^{m}(\nabla \Psi _{\ast
})(e_{i},e_{i}),  \label{EQ2.18}
\end{equation}%
where $\{e_{1},e_{2},......,e_{m}\}$ is the orthonormal basis on $M^{\ast }$%
. Then it follows that ns condition for $\Psi $ to be harmonic is $\tau
(\Psi )=0$ \cite{Baird}.

Now, we recall an important Lemma from \cite{Baird}, which will be needed in
the study of whole paper.

\begin{lem}
\label{LE3} Let $\Psi :(M^{\ast},g_{m})\rightarrow (N^{\ast},g_{n})$ be a
horizontal conformal submersion, $P,Q$ are horizontal vector fields and $V,W$
are vertical vector fields. Then
\end{lem}

$(a)$ $(\nabla \Psi _{\ast })(P,Q)=P(\ln \lambda_{\ast} )\Psi _{\ast
}(Q)+Q(\ln \lambda_{\ast} )\Psi _{\ast }(P)-g_{m}(P,Q)\Psi _{\ast }({grad}\
ln\lambda_{\ast} ),$

$(b)$ $(\nabla \Psi _{\ast })(V,W)=-\Psi _{\ast }({\mathcal{T}}_{V}W),$

$(c)$ $(\nabla \Psi _{\ast })(P,V)=-\Psi _{\ast }(\nabla _{P}^{M}V)=-\Psi
_{\ast }(A_{P}V).$

\section{ Semi-invariant conformal $\protect\zeta ^{\perp }$-Riemannian
submersions \label{sect-con}}

In this segment, we deal with the definition and examples of sic $\zeta
^{\perp }$-Rs from acm manifolds onto Riemannian manifolds. We acquire the
integrability of distributions and also geometry of leaves of $ker\Psi
_{\ast }$ and $(ker\Psi _{\ast })^{\perp }$. \newline

\begin{defn}
\label{DEF3.1}Let $\left( M^{\ast},\varphi_{\ast} ,\zeta ,\eta_{\ast}
,g_{m}\right) $ be an acm manifold and $\left( N^{\ast},g_{n}\right) $ be a
Riemannian manifold. A horizontally conformal submersion $\Psi
:(M^{\ast},g_{m})\rightarrow (N^{\ast},g_{n})$ with slope $\lambda_{\ast} $,
is said to be a sic $\zeta ^{\perp }$-Rs if $\zeta $ is normal to $ker\Psi
_{\ast }$ and there is a distribution $D^{\ast}_{1}\subseteq ker\Psi _{{\ast
}}$ such that
\begin{equation}
ker\Psi _{\ast }=D^{\ast}_{1}\oplus D^{\ast}_{2},  \label{EQ3.1}
\end{equation}%
\begin{equation}
\varphi_{\ast} (D^{\ast}_{1})=D^{\ast}_{1},\ \varphi_{\ast}
(D^{\ast}_{2})\subseteq (ker\Psi _{\ast })^{\perp },  \label{EQ3.2}
\end{equation}%
where $D^{\ast}_{2}$ is orthogonal complementary to $D^{\ast}_{1}$ in $%
ker\Psi _{{\ast }}$.
\end{defn}

It is to be noted that the distribution $ker\Psi _{\ast }$ is integrable.%
\newline

Further to prove the consistency of sic $\zeta ^{\perp }$-Rs in acm
manifolds, we are putting some examples here.

\begin{ex}
Every anti-invariant $\zeta ^{\perp }$-Riemannian submersion from an acm
manifold onto a Riemannian manifold is a sic $\zeta ^{\perp }$-Rs with $%
\lambda_{\ast} =I$ and $D^{\ast}_{1}=\{0\}$, where $I$ denotes the identity
function.
\end{ex}

\begin{ex}
Every semi-invariant $\zeta ^{\perp }$ submersion from an acm manifold onto
a Riemannian manifold is a sic $\zeta ^{\perp }$-Rs with $\lambda_{\ast} =I.$
\end{ex}

\begin{ex}
\cite{Blair} Let $(\mathbb{R}^{m},g,\varphi _{\ast },\zeta ,\eta _{\ast })$,
$(m=2n+1)$ be a Sasakian manifold given by%
\[
\eta _{\ast }=\frac{1}{2}\left( dw-\sum\limits_{i=1}^{n}v^{i}du^{i}\right)
,\quad \zeta =2\frac{\partial }{\partial w},
\]%
\[
g=\frac{1}{4}\eta _{\ast }\otimes \eta _{\ast
}+\sum\limits_{i=1}^{n}(du^{i}\otimes du^{i}+dv^{i}\otimes dv^{i}),
\]%
\[
\varphi _{\ast }=\left(
\begin{array}{ccc}
0 & \delta _{ij} & 0 \\
-\delta _{ij} & 0 & 0 \\
0 & v^{j} & 0%
\end{array}%
\right) ,
\]%
where $\{(u^{i},v^{i},w)|i=1,\ldots ,n\}$ are the cartesian coordinates. The
vector fields $X_{i}=2\frac{\partial }{\partial v^{i}},\quad
X_{n+i}=2\left( \frac{\partial }{\partial u^{i}}+v^{i}\frac{\partial }{%
\partial w}\right) $ and $\zeta $ form a $\varphi _{\ast }$-basis for the
contact metric structure.

\noindent Let $\Psi :\left( \mathbb{R}^{7},g_{1}\right) \rightarrow \left(
\mathbb{R}^{2},g_{2}\right) $ be a Riemannian submersion defined by%
\[
\Psi (u_{1},u_{2},u_{3},v_{1},v_{2},v_{3},w)=(\sinh v_{3}\cos w,\cosh
v_{3}\sin w).
\]%
Then it follows that
\[
ker\Psi _{\ast }=span\{L_{1}=\partial u_{1},L_{2}=\partial
u_{2},L_{3}=\partial u_{3},L_{4}=\partial v_{1},L_{5}=\partial v_{2}\}
\]%
and%
\[
(ker\Psi _{\ast })^{\bot }=span\left\{ W_{1}=\partial v_{3},W_{2}=\partial
w\right\} .
\]%
Hence, we have $\varphi _{\ast }L_{1}=-L_{4},\varphi _{\ast
}L_{2}=-L_{5},\varphi _{\ast }L_{3}=-W_{1}$, $\varphi _{\ast }L_{4}=L_{1}$
and $\varphi _{\ast }L_{5}=L_{2}.$ Thus it follows that $D_{1}^{\ast
}=span\{L_{1},L_{2},L_{4},L_{5}\}$, $D_{2}^{\ast }=span\{L_{3}\}$ and $\zeta
\in (ker\Psi _{\ast })^{\bot }$. We can easily compute that
\begin{eqnarray*}
g_{1}\left( \Psi _{\ast }W_{1},\Psi _{\ast }W_{1}\right)  &=&\left( \sinh
^{2}v_{3}\sin ^{2}w+\cosh ^{2}v_{3}\cos ^{2}w\right) g_{2}\left(
W_{1},W_{1}\right) , \\
g_{1}\left( \Psi _{\ast }W_{2},\Psi _{\ast }W_{2}\right)  &=&\left( \sinh
^{2}v_{3}\sin ^{2}w+\cosh ^{2}v_{3}\cos ^{2}w\right) g_{2}\left(
W_{2},W_{2}\right) .
\end{eqnarray*}%
Therefore, $\Psi $ is a sic $\zeta ^{\perp }$-Rs with $\lambda _{\ast }=%
\sqrt{\left( \sinh ^{2}v_{3}\sin ^{2}w+\cosh ^{2}v_{3}\cos ^{2}w\right) }.$
\end{ex}

Let $\Psi $ be a sic $\zeta ^{\perp }$-Rs from acm manifolds $%
(M^{\ast},\varphi_{\ast} ,\zeta ,\eta_{\ast} ,g_{m})$ onto Riemannian
manifolds $(N^{\ast},g_{n})$. Consider
\[
\left( \ker \Psi _{\ast }\right) ^{\bot }=\varphi_{\ast} D^{\ast}_{2}\oplus
\mu ,
\]%
where $\mu $ is the complementary distribution to $\varphi_{\ast}
D^{\ast}_{2}$ in $(ker\Psi _{\ast })^{\bot }$ and $\varphi_{\ast} \mu
\subset \mu $. Therefore $\zeta \in \mu $.

For $V\in \Gamma (ker\Psi _{\ast })$, we write
\begin{equation}
\varphi_{\ast} V=\phi V+\omega V,  \label{EQ3.3}
\end{equation}%
where $\phi V\in \Gamma (D^{\ast}_{1})$ and $\omega V\in \Gamma
(\varphi_{\ast} D^{\ast}_{2})$.

For $P\in \Gamma(ker\Psi _{\ast })^{\perp }$, we have
\begin{equation}
\varphi_{\ast} P=\alpha P+\beta P,  \label{EQ3.4}
\end{equation}%
where $\alpha P\in \Gamma (D^{\ast}_{2})$ and $\beta P\in \Gamma (\mu )$.

Now, using (\ref{EQ2.10}), (\ref{EQ2.11}), (\ref{EQ3.3}) and (\ref{EQ3.4}),
we get
\begin{equation}
(\nabla _{V}^{M}\phi )W=\alpha T_{V}W-T_{V}\omega W,  \label{EQ3.5}
\end{equation}%
\begin{equation}
(\nabla _{V}^{M}\omega )W=\beta T_{V}W-T_{V}\phi W  \label{EQ3.6}
\end{equation}%
for $V,W\in \Gamma (ker\Psi _{\ast })$, where
\begin{equation}
(\nabla _{V}^{M}\phi )W=\hat{\nabla}_{V}\phi W-\phi \hat{\nabla}_{V}W
\label{EQ3.7}
\end{equation}%
and
\begin{equation}
(\nabla _{V}^{M}\omega )W={\mathcal{H}_{1}}\nabla _{V}^{M}\omega W-\omega
\hat{\nabla}_{V}W.  \label{EQ3.8}
\end{equation}
Now, we investigate the integrability of the distributions $D^{\ast}_{1}$
and $D^{\ast}_{2} $.

\begin{lem}
Let $\Psi $ be sic $\zeta ^{\perp }$-Rs from Sasakian manifolds $%
(M^{\ast},\varphi_{\ast} ,\zeta ,\eta_{\ast} ,g_{m})$ onto Riemannian
manifolds $(N^{\ast},g_{n})$. Then the ns condition for the distribution $%
D^{\ast}_{1}$ to be integrable is
\[
(\nabla \Psi _{\ast })(V,\varphi_{\ast} U)-(\nabla \Psi _{\ast
})(U,\varphi_{\ast} V)\in \Gamma (\Psi _{\ast }(\mu ))
\]%
for $U,V\in \Gamma (D^{\ast}_{1})$.
\end{lem}

\noindent {\textbf Proof}. Let $U,V\in \Gamma (D_{1}^{\ast })$ and $W\in
\Gamma (D_{2}^{\ast })$. With the help of (\ref{EQ2.2}), we get
\begin{equation}
g_{m}(\nabla _{U}\varphi _{\ast }V,\varphi _{\ast }W)=g_{m}(\nabla _{U}V,W).
\label{Eq 3.9}
\end{equation}%
So
\[
g_{m}([U,V],W)=g_{m}(\nabla _{U}V,W)-g_{m}(\nabla _{V}U,W).
\]%
Using (\ref{EQ2.10}) and (\ref{Eq 3.9}), we obtain
\begin{eqnarray}
g_{m}([U,V],W) &=&g_{m}(\nabla _{U}\varphi _{\ast }V,\varphi _{\ast
}W)-g_{m}(\nabla _{V}\varphi _{\ast }U,\varphi _{\ast }W)  \nonumber \\
&=&g_{m}(T_{U}\varphi _{\ast }V-T_{V}\varphi _{\ast }U,\varphi _{\ast }W).
\nonumber
\end{eqnarray}%
The distribution $D_{1}^{\ast }$ is integrable if and only if $%
g_{m}([U,V],W)=g_{m}([U,V],Z)=0$ for $U,V\in \Gamma (D_{1}^{\ast })$, $W\in
\Gamma (D_{2}^{\ast })$ and $Z\in \Gamma ((ker\Psi _{\ast })^{\perp })$.
Since $ker\Psi _{\ast }$ is integrable, so we immediately have $%
g_{m}([U,V],Z)=0$. Thus the ns condition for the distribution $D_{1}^{\ast }$
to be integrable is $g_{m}([U,V],W)=0$. Since $\Psi $ is a conformal
submersion, in view of Lemma \ref{LE3} $(b)$, we conclude that
\[
g_{m}([U,V],W)=\frac{1}{\lambda _{\ast }^{2}}g_{n}((\nabla \Psi _{\ast
})(V,\varphi _{\ast }U)-(\nabla \Psi _{\ast })(U,\varphi _{\ast }V),\Psi
_{\ast }\varphi _{\ast }W).
\]

\begin{lem}
Let $\Psi $ be sic $\zeta ^{\perp }$-Rs from Sasakian manifolds $%
(M^{\ast},\varphi_{\ast} ,\zeta ,\eta_{\ast} ,g_{m})$ onto Riemannian
manifolds $(N^{\ast},g_{n})$. Then the ns condition for the integrability of
distribution $D^{\ast}_{2}$ is
\[
\hat{\nabla}_{W}\phi \xi+T_{W}\omega \xi-\hat{\nabla}_{\xi}\phi
W-T_{\xi}\omega W\in \Gamma (D^{\ast}_{2})
\]%
for $W,\xi \in \Gamma (D^{\ast}_{2})$.
\end{lem}

\noindent {\textbf Proof}. Let $U\in \Gamma (D^{\ast}_{1})$ and $W,\xi\in
\Gamma (D^{\ast}_{2})$. Using (\ref{EQ2.2}), (\ref{EQ2.3}), (\ref{EQ2.9}), (%
\ref{EQ2.10}), (\ref{EQ2.11}) and (\ref{Eq 3.9}), we get

\begin{eqnarray*}
g_{m}(\varphi_{\ast} \lbrack W,\xi],U) &=&g_{m}(\nabla _{W}\varphi_{\ast}
\xi,U)-g_{m}(\nabla _{\xi}\varphi_{\ast} W,U) \\
&=&g_{m}(T_{W}\phi \xi+\hat{\nabla}_{W}\phi \xi+T_{W}\omega \xi+{\mathcal{H}%
_{1}}\nabla _{W}\omega \xi \\
&&-T_{\xi}\phi W-\hat{\nabla}_{\xi}\phi W-T_{\xi}\omega W-{\mathcal{H}_{1}}%
\nabla _{\xi}\omega W,U) \\
&=&g_{m}(\hat{\nabla}_{W}\phi \xi+T_{W}\omega \xi-\hat{\nabla}_{\xi}\phi
W-T_{\xi}\omega W,U).
\end{eqnarray*}

\begin{th}
Let $\Psi $ be a sic $\zeta ^{\perp }$-Rs from Sasakian manifolds $%
(M^{\ast},\varphi_{\ast} ,\zeta ,\eta_{\ast} ,g_{m})$ onto Riemannian
manifolds $(N^{\ast},g_{n})$. Then the ns conditions for the distribution $%
(ker\Psi _{\ast })^{\perp }$ to be integrable are
\[
A_{Q}\omega \alpha P-A_{P}\omega \alpha Q-\varphi_{\ast} A_{P}\beta
Q+\varphi_{\ast} A_{Q}\beta P\notin \Gamma (D^{\ast}_{1})
\]%
and
\begin{eqnarray*}
\frac{1}{\lambda_{\ast} ^{2}}g_{n}(\nabla _{V}\Psi _{\ast }CU-\nabla
_{U}\Psi _{\ast }CV,\Psi _{\ast }\varphi_{\ast} Z)&
=g_{m}(A_{V}BU-A_{U}BV-CV(ln\lambda_{\ast} )U \\
& +CU(ln\lambda_{\ast} )V+2g_{m}(U,CV)grad\ ln\lambda_{\ast} \\
& +\eta_{\ast} (U)V-\eta_{\ast} (V)U,\varphi_{\ast} Z)
\end{eqnarray*}%
for $U,V\in \Gamma (ker\Psi _{\ast })^{\perp })$ and $Z\in \Gamma
(D^{\ast}_{2})$.
\end{th}

\noindent {\textbf Proof.} The ns condition for the distribution $(ker\Psi
_{\ast })^{\perp }$ to be integrable on manifold $M^{\ast}$ is that
\[
g_{m}([P,Q],U)=0\ \ {\rm and}\ \ g_{m}([P,Q],W)=0,
\]%
for $P,Q\in \Gamma ((ker\Psi _{\ast })^{\perp })$, $U\in \Gamma
(D^{\ast}_{1})$ and $W\in \Gamma (D^{\ast}_{2})$.

Using (\ref{EQ2.1}), (\ref{EQ2.2}), (\ref{EQ2.4}) and (\ref{EQ3.4}), we
infer
\begin{eqnarray*}
g_{m}([P,Q],U)& =g_{m}(\varphi_{\ast} \lbrack P,Q],\varphi_{\ast} U) \\
& =g_{m}(\varphi_{\ast} \nabla _{P}Q-\varphi_{\ast} \nabla
_{Q}P,\varphi_{\ast} U) \\
& =g_{m}(\nabla _{P}\alpha Q,\varphi_{\ast} U)+g_{m}(\nabla _{P}\beta
Q,\varphi_{\ast} U) \\
& -g_{m}(\nabla _{Q}\alpha P,\varphi_{\ast} U)-g_{m}(\nabla _{Q}\beta
P,\varphi_{\ast} U) \\
& =g_{m}(\varphi_{\ast} \nabla _{P}\alpha Q,U)-g_{m}(\varphi_{\ast} \nabla
_{P}\beta Q,U) \\
& +g_{m}(\varphi_{\ast} \nabla _{Q}\alpha P,U)+g_{m}(\varphi_{\ast} \nabla
_{Q}\beta P,U).
\end{eqnarray*}%
By simple calculations, we have
\begin{eqnarray*}
g_{m}([P,Q],U)& =-g_{m}(\nabla _{P}\varphi_{\ast} \alpha Q,U)+g_{m}(\nabla
_{P}\beta Q,\varphi_{\ast} U) \\
& +g_{m}(\nabla _{Q}\varphi_{\ast} \alpha P,U)-g_{m}(\nabla _{Q}\beta
P,\varphi_{\ast} U).
\end{eqnarray*}%
By (\ref{EQ3.3}), we have $\varphi_{\ast} \alpha P=\phi \alpha P+\omega
\alpha P$. Since $\varphi_{\ast} \alpha P$, $\omega \alpha P\in \Gamma
((ker\Psi _{\ast })^{\perp })$, so $\phi \alpha P=0$, $\forall P\in \Gamma
((ker\Psi _{\ast })^{\perp })$. Using (\ref{EQ2.4}),  (\ref{EQ2.13}) and (%
\ref{EQ3.3}), we get
\begin{equation}
g_{m}([P,Q],U)=g_{m}(A_{Q}\omega \alpha P-A_{P}\omega \alpha
Q-\varphi_{\ast} A_{P}\beta Q+\varphi_{\ast} A_{Q}\beta P,U).  \label{EQ3.9}
\end{equation}%
With the help of (\ref{EQ2.1}), (\ref{EQ2.2}), (\ref{EQ2.4}) and (\ref{EQ3.4}%
), we deduce
\begin{eqnarray*}
g_{m}([P,Q],W)& =g_{m}(\nabla _{P}\alpha Q,\varphi_{\ast} W)+g_{m}(\nabla
_{P}\beta Q,\varphi_{\ast} W) \\
& -g_{m}(\nabla _{Q}\alpha P,\varphi_{\ast} W)-g_{m}(\nabla _{Q}\beta
P,\varphi_{\ast} W) \\
& +\eta_{\ast} (P)g_{m}(Q,\varphi_{\ast} W)-\eta_{\ast}
(Q)g_{m}(P,\varphi_{\ast} W)
\end{eqnarray*}%
Since $\Psi $ is sic $\zeta ^{\perp }$-Rs, using (\ref{EQ2.12}), (\ref%
{EQ2.13}), (\ref{EQ2.15}) and Lemma \ref{LE2}, \ref{LE3}, we obtain
\begin{eqnarray*}
g_{m}([P,Q],W)& =-\eta_{\ast} (Q)g_{m}(P,\varphi_{\ast} W)+\eta_{\ast}
(P)g_{m}(Q,\varphi_{\ast} W) \\
& -\frac{1}{\lambda_{\ast} ^{2}}g_{n}((\nabla \Psi _{\ast })(P,\alpha
Q),\Psi _{\ast }\varphi_{\ast} W)+\frac{1}{\lambda_{\ast} ^{2}}g_{n}((\nabla
\Psi _{\ast })(Q,\alpha P),\Psi _{\ast }\varphi_{\ast} W) \\
& +\frac{1}{\lambda_{\ast} ^{2}}g_{n}\{-P(ln\lambda_{\ast} )\Psi _{\ast
}\beta Q-\beta Q(ln\lambda_{\ast} )\Psi _{\ast }P \\
& +g_{m}(P,\beta Q)\Psi _{\ast }(grad\ ln\lambda_{\ast} )+\nabla _{P}^{\Psi
}\Psi _{\ast }\beta Q,\Psi _{\ast }\varphi_{\ast} W\} \\
& -\frac{1}{\lambda_{\ast} ^{2}}g_{n}\{-Q(ln\lambda_{\ast} )\Psi _{\ast
}\beta P-\beta P(ln\lambda_{\ast} )\Psi _{\ast }Q \\
& +g_{m}(Q,\beta P)\Psi _{\ast }(grad\ ln\lambda_{\ast} )+\nabla _{Q}^{\Psi
}\Psi _{\ast }\beta P,\Psi _{\ast }\varphi_{\ast} W\},
\end{eqnarray*}%
\begin{eqnarray*}
g_{m}([P,Q],W)& =-\eta_{\ast} (Q)g_{m}(P,\varphi_{\ast} W)+\eta_{\ast}
(P)g_{m}(Q,\varphi_{\ast} W)-g_{m}(A_{P}\alpha Q,\varphi_{\ast} W) \\
&+g_{m}(A_{Q}\alpha P,\varphi_{\ast} W) -\frac{1}{\lambda_{\ast} ^{2}}%
\{g_{m}(grad\ ln\lambda_{\ast} ,P)g_{n}(\Psi _{\ast }\beta Q,\Psi _{\ast
}\varphi_{\ast} W) \\
&+g_{m}(grad\ ln\lambda_{\ast} ,\beta Q)g_{n}(\Psi _{\ast }P,\Psi _{\ast
}\varphi_{\ast} W) \\
&-g_{m}(P,\beta Q)g_{n}(\Psi _{\ast }(grad\ ln\lambda_{\ast} ),\Psi _{\ast
}\varphi_{\ast} W) \\
&-g_{m}(grad\ ln\lambda_{\ast} ,Q)g_{n}(\Psi _{\ast }\beta P,\Psi _{\ast
}\varphi_{\ast} W) \\
&-g_{m}(grad\ ln\lambda_{\ast} ,\beta P)g_{n}(\Psi _{\ast }Q,\Psi _{\ast
}\varphi_{\ast} W) \\
&+g_{m}(Q,\beta P)g_{n}(\Psi _{\ast }(grad\ ln\lambda_{\ast} ),\Psi _{\ast
}\varphi_{\ast} W) \\
&-g_{n}(\nabla _{P}^{\Psi }\Psi _{\ast }\beta Q,\Psi _{\ast }\varphi_{\ast}
W)+g_{n}(\nabla _{Q}^{\Psi }\Psi _{\ast }\beta P,\Psi _{\ast }\varphi_{\ast}
W)\}.
\end{eqnarray*}%
Using definition \ref{DEF3.1}, we have
\begin{eqnarray}
g_{m}([P,Q],W) &=&-\eta_{\ast} (Q)g_{m}(P,\varphi_{\ast} W)+\eta_{\ast}
(P)g_{m}(Q,\varphi_{\ast} W)  \nonumber \\
&&+g_{m}(A_{Q}\alpha P-A_{P}\alpha Q-\left( \beta Q\right) (P\ln
\lambda_{\ast} )  \nonumber \\
&&+\left( \beta P\right) (Q\ln \lambda_{\ast} )+2g_{m}(P,\beta Q)grad\
ln\lambda_{\ast} ,\varphi_{\ast} W)  \nonumber \\
&&+g_{m}((\beta Pln\lambda_{\ast} )Q-(\beta Qln\lambda_{\ast}
)P,\varphi_{\ast} W)  \nonumber \\
&&+\frac{1}{\lambda_{\ast} ^{2}}g_{n}(\nabla _{P}^{\Psi }\Psi _{\ast }\beta
Q-\nabla _{Q}^{\Psi }\Psi _{\ast }\beta P,\Psi _{\ast }\varphi_{\ast} W).
\label{EQ3.10.}
\end{eqnarray}%
From the above expression, we entail that
\begin{eqnarray}
\frac{1}{\lambda_{\ast} ^{2}}g_{n}(\nabla _{Q}^{\Psi }\Psi _{\ast }\beta
P-\nabla _{P}^{\Psi }\Psi _{\ast }\beta Q,\Psi _{\ast }\varphi_{\ast} W)&
=g_{m}(A_{Q}\alpha P-A_{P}\alpha Q-\left( \beta Q\right) (P\ln
\lambda_{\ast} )  \nonumber  \label{EQ3.10} \\
& +\left( \beta P\right) (Q\ln \lambda_{\ast} )+2g_{m}(P,\beta Q)grad\
ln\lambda_{\ast}  \nonumber \\
& +(\beta Pln\lambda_{\ast} )Q-(\beta Qln\lambda_{\ast} )P+\eta_{\ast} (P)Q
\nonumber \\
&-\eta_{\ast} (Q)P,\varphi_{\ast} W).
\end{eqnarray}%
So from (\ref{EQ3.9}) and (\ref{EQ3.10}), we obtain the results.\newline
Now, we find out the required necessary and sufficient conditions for sic $%
\zeta ^{\perp }$-Rs to be homothetic map.

\begin{th}
Let $\Psi $ be a sic $\zeta ^{\perp }$-Rs from Sasakian manifolds $%
(M^{\ast},\varphi_{\ast} ,\zeta ,\eta_{\ast} ,g_{m})$ onto Riemannian
manifolds $(N^{\ast},g_{n})$ with integrable distribution $(ker\Psi _{\ast
})^{\perp }$. Then the ns condition for the $\Psi $ to be homothetic map is
\begin{eqnarray}
g_{n}(\nabla _{Q}^{\Psi }\Psi _{\ast }\beta P-\nabla _{P}^{\Psi }\Psi _{\ast
}\beta Q,\Psi _{\ast }\varphi_{\ast} W)& =\lambda_{\ast}
^{2}g_{m}(A_{Q}\alpha P-A_{P}\alpha Q  \nonumber \\
& +\eta_{\ast} (P)Q-\eta_{\ast} (Q)P,\varphi_{\ast} W)  \label{EQ3.11}
\end{eqnarray}%
for $P,Q\in \Gamma ((ker\Psi _{\ast })^{\perp })$ and $W\in \Gamma
(D^{\ast}_{2})$.
\end{th}

\noindent {\textbf Proof.} Let $P,Q\in \Gamma ((ker\Psi _{\ast })^{\perp })$
and $W\in \Gamma (D^{\ast}_{2})$. Using (\ref{EQ3.10}), we infer
\begin{eqnarray}
g_{m}([P,Q],W)& =-\eta_{\ast} (Q)g_{m}(P,\varphi_{\ast} W)+\eta_{\ast}
(P)g_{m}(Q,\varphi_{\ast} W)  \nonumber \\
& +g_{m}(A_{Q}\alpha P-A_{P}\alpha Q-\left( \beta Q\right) (P\ln
\lambda_{\ast} )  \nonumber \\
& +\left( \beta P\right) (Q\ln \lambda_{\ast} )+2g_{m}(P,\beta Q)grad\
ln\lambda_{\ast} ,\varphi_{\ast} W)  \nonumber \\
& +g_{m}((\beta Pln\lambda_{\ast} )Q-(\beta Qln\lambda_{\ast}
)P,\varphi_{\ast} W)  \nonumber \\
& +\frac{1}{\lambda_{\ast} ^{2}}g_{n}(\nabla _{P}^{\Psi }\Psi _{\ast }\beta
Q-\nabla _{Q}^{\Psi }\Psi _{\ast }\beta P,\Psi _{\ast }\varphi_{\ast} W).
\end{eqnarray}%
If $\Psi $ is a parallel homothetic map, therefore we acquire (\ref{EQ3.11}%
). Conversely, if (\ref{EQ3.11}) holds, then we gain
\begin{eqnarray}
0& =g_{m}((\beta Pln\lambda_{\ast} )Q-(\beta Qln\lambda_{\ast}
)P,\varphi_{\ast} W)  \nonumber \\
& +g_{m}(-\left( \beta Q\right) (P\ln \lambda_{\ast} )+\left( \beta P\right)
(Q\ln \lambda_{\ast} )+2g_{m}(P,\beta Q)grad\ ln\lambda_{\ast}
,\varphi_{\ast} W).  \label{EQ3.12}
\end{eqnarray}%
Replacing $Q$ by $\varphi_{\ast} W$ for $W\in \Gamma (D^{\ast}_{2})$ in (\ref%
{EQ3.12}), we entail that
\[
g_{m}({grad}\ln \lambda_{\ast} ,\beta P)g_{m}(\varphi_{\ast}
W,\varphi_{\ast} W)=0,
\]%
which implies that $\lambda_{\ast} $ is a constant on $\Gamma (\mu )$ since $%
||\varphi_{\ast} W||\neq 0$. \newline
Now, again replace $Q$ by $\beta P$ for $P\in \Gamma (\mu )$ in (\ref{EQ3.12}%
), we have
\[
g_{m}(grad\ ln\lambda_{\ast} ,\varphi_{\ast} W)g_{m}(P,\beta
^{2}P)=g_{m}(grad\ ln\lambda_{\ast} ,\varphi_{\ast} W)g_{m}(\beta P,\beta
P)=0.
\]%
Now, we see that the $||\beta P||\neq 0$, this concludes $\lambda_{\ast} $
is a constant on $\Gamma (\varphi_{\ast} D^{\ast}_{2})$, which finalize the
proof of our theorem.

Now, we come to the following results for the geometry of leaves of
horizontal distributions.

\begin{th}
Let $\Psi $ be a sic $\zeta ^{\perp }$-Rs from Sasakian manifolds $%
(M^{\ast},\varphi_{\ast} ,\zeta ,\eta_{\ast} ,g_{m})$ onto Riemannian
manifolds $(N^{\ast},g_{n})$. Then the ns conditions for the $(ker\Psi
_{\ast }^{\perp })$ to be completely geodesic foliation on $M^{\ast}$ are
\[
A_{P}\varphi_{\ast} Q+{\mathcal{V}_{1}}\nabla _{P}^{m}\alpha Q\in \Gamma
(D^{\ast}_{2}),
\]%
and
\begin{eqnarray*}
g_{n}(\nabla _{P}\Psi _{\ast }\varphi_{\ast} W,\Psi _{\ast }\beta Q)=&
\lambda_{\ast} ^{2}g_{m}(-A_{P}\alpha Q+\beta Q(ln\lambda_{\ast}
)P-g_{m}(P,\beta Q)({grad}\ln \lambda_{\ast} ) \\
& -P(ln\lambda_{\ast} )\beta Q+\eta_{\ast} (Q)P,\varphi_{\ast} W)
\end{eqnarray*}%
for $P,Q\in \Gamma (ker\Psi _{\ast })^{\perp }$, $U\in \Gamma (D^{\ast}_{1})$
and $W\in \Gamma (D^{\ast}_{2})$.
\end{th}

\noindent {\textbf Proof.} Here, it is to be noted that the distribution $%
(ker\Psi _{\ast })^{\perp }$ characterizes a completely geodesic foliation
on $M$ if and only if $g_{m}(\nabla _{P}^{m}Q,U)=0$ and $g_{m}(\nabla
_{P}^{m}Q,W)=0$ for $P,Q\in \Gamma (ker\Psi _{\ast })^{\perp }$, $U\in
\Gamma (D^{\ast}_{1})$ and $W\in \Gamma (D^{\ast}_{2})$. Then with the help
of (\ref{EQ2.1}),(\ref{EQ2.2}), (\ref{EQ2.4}), (\ref{EQ2.12}), (\ref{EQ2.13}%
) and (\ref{EQ3.4}), we infer
\begin{eqnarray}
g_{m}(\nabla _{P}^{m}Q,U)& =g_{m}(\nabla _{P}^{m}\varphi_{\ast}
Q,\varphi_{\ast} U)-\eta_{\ast} (Q)g_{m}(P,\varphi_{\ast} U)  \label{EQ3.13}
\\
& =-g_{m}(\varphi_{\ast} \left( A_{P}\varphi_{\ast} Q+{\mathcal{V}_{1}}%
\nabla _{P}^{m}\alpha Q\right) ,U),  \nonumber
\end{eqnarray}%
which implies that $A_{P}\varphi_{\ast} Q+{\mathcal{V}_{1}}\nabla
_{P}^{m}\alpha Q\in \Gamma (D^{\ast}_{2})$.\newline
Using (\ref{EQ2.1}),(\ref{EQ2.2}), (\ref{EQ2.4}) and (\ref{EQ3.4}), we
observe that
\begin{eqnarray*}
g_{m}(\nabla _{P}^{m}Q,W)& =g_{m}(\nabla _{P}^{m}\alpha Q+\nabla
_{P}^{m}\beta Q,\varphi_{\ast} W) \\
& -\eta_{\ast} (Q)g_{m}(P,\varphi_{\ast} W), \\
& =-g_{m}(\alpha Q,\nabla _{P}^{m}\varphi_{\ast} W)-g_{m}(\beta Q,\nabla
_{P}^{m}\varphi_{\ast} W) \\
& -\eta_{\ast} (Q)g_{m}(P,\varphi_{\ast} W).
\end{eqnarray*}%
Since $\Psi $ is a conformal Riemannian submersion, in view of (\ref{EQ2.13}%
), (\ref{EQ2.15}) and Lemma \ref{LE2}, we gain
\begin{eqnarray*}
g_{m}(\nabla _{P}^{m}Q,U)& =-g_{m}(\alpha Q,A_{P}\varphi_{\ast} W)+\frac{1}{%
\lambda_{\ast} ^{2}}g_{m}({grad}\ln \lambda_{\ast} ,\varphi_{\ast}
W)g_{n}(\Psi _{\ast }P,\Psi _{\ast }\beta Q) \\
& +\frac{1}{\lambda_{\ast} ^{2}}g_{m}({grad}\ln \lambda_{\ast} ,P)g_{n}(\Psi
_{\ast }\varphi_{\ast} W,\Psi _{\ast }\beta Q,) \\
&-\frac{1}{\lambda_{\ast} ^{2}}g_{n}(\Psi _{\ast }({grad}\ln \lambda_{\ast}
),\Psi _{\ast }\beta Q)g_{m}(P,\varphi_{\ast} W) \\
& +\frac{1}{\lambda_{\ast} ^{2}}g_{n}(\nabla _{P}^{\Psi }\Psi _{\ast
}\varphi_{\ast} W,\Psi _{\ast }\beta Q)-\eta_{\ast}
(Q)g_{m}(P,\varphi_{\ast} W).
\end{eqnarray*}%
Using the definition of sic $\zeta ^{\perp }$-Rs, we obtain
\begin{eqnarray}
g_{m}(\nabla _{P}^{m}Q,U)& =g_{m}(A_{P}\alpha Q-\beta Q(ln\lambda_{\ast}
)P+g_{m}(P,\beta Q)({grad}\ln \lambda_{\ast} )  \label{EQ3.14} \\
& +P(ln\lambda_{\ast} )\beta Q-\eta_{\ast} (Q)P,\varphi_{\ast} W)+\frac{1}{%
\lambda_{\ast} ^{2}}g_{n}(\nabla _{P}\Psi _{\ast }\varphi_{\ast} W,\Psi
_{\ast }\beta Q).  \nonumber
\end{eqnarray}%
Thus, we obtain the proof from (\ref{EQ3.13}) and (\ref{EQ3.14}).

Now, we state the following definition for further result.

\begin{defn}
Let $\Psi $ be a sic $\zeta ^{\perp }$-Rs from acm manifolds $(M^{\ast
},\varphi _{\ast },\zeta ,\eta _{\ast },g_{m})$ onto Riemannian manifolds $%
(N^{\ast },g_{n})$. The distribution $D_{2}^{\ast }$ is parallel along $%
(ker\Psi _{\ast })^{\perp }$ if $\nabla _{P}^{m}W\in \Gamma (D_{2}^{\ast })$
for $P\in \Gamma ((ker\Psi _{\ast })^{\perp })$ and $W\in \Gamma
(D_{2}^{\ast })$.
\end{defn}

\begin{cor}
Let $\Psi $ be a sic $\zeta ^{\perp }$-Rs from Sasakian manifolds $%
(M^{\ast},\varphi_{\ast} ,\zeta ,\eta_{\ast} ,g_{m})$ onto Riemannian
manifolds $(N^{\ast},g_{n})$ such that $D^{\ast}_{2}$ is parallel along $%
(ker\Psi _{\ast })^{\perp }$. Then the ns condition for the $\Psi $ to be
horizontally homothetic map is
\begin{equation}
\lambda_{\ast} ^{2}g_{m}(A_{P}\alpha Q,\varphi_{\ast} W)=g_{n}(\nabla
_{P}\Psi _{\ast }\varphi_{\ast} W,\Psi _{\ast }\beta Q)  \label{EQ3.15}
\end{equation}%
for $P,Q\in \Gamma (ker\Psi _{\ast })^{\perp }$ and $W\in \Gamma
(D^{\ast}_{2})$.
\end{cor}

\noindent {\textbf Proof.} Using (\ref{EQ2.13}), (\ref{EQ3.15}) and Lemma \ref%
{LE3}, we obtain
\begin{equation}
-g_{m}(grad\ ln\lambda_{\ast} ,\beta Q)g_{m}(P,\varphi_{\ast}
W)+g_{m}(P,\beta Q)g_{m}(grad\ ln\lambda_{\ast} ,\varphi_{\ast} W)=0.
\label{EQ3.16}
\end{equation}%
Now, replacing $P$ by $\varphi_{\ast} W$ for $W\in \Gamma (D^{\ast}_{2})$ in %
\ref{EQ3.16}, we get
\[
g_{m}(grad\ ln\lambda_{\ast} ,\beta Q)g_{m}(\varphi_{\ast} W,\varphi_{\ast}
W)=0.
\]%
Since $||\varphi_{\ast} W||\neq 0$, Thus $\lambda_{\ast} $ is a constant on $%
\Gamma (\mu )$.\newline
Now again replace $P$ by $\beta Q$ for $Q\in \Gamma (ker\Psi _{\ast
})^{\perp }$ in (\ref{EQ3.16}), we get
\[
g_{m}(\beta Q,\beta Q)g_{m}(grad\ ln\lambda_{\ast} ,\varphi_{\ast} W)=0.
\]%
Since $||\beta Q||\neq 0$, Thus $\lambda_{\ast} $ is a constant on $\Gamma
(\varphi_{\ast} D^{\ast}_{2})$. Therefore $\lambda_{\ast} $ is constant on $%
\Gamma (ker\Psi _{\ast })^{\perp }$.\newline
The converse easily follows from (\ref{EQ3.16}).\newline

Now, we explore the geometric investigation of leaves of the distribution $%
(ker\Psi _{\ast }) $.

\begin{th}
Let $\Psi $ be a sic $\zeta ^{\perp }$-Rs from Sasakian manifolds $%
(M^{\ast},\varphi_{\ast} ,\zeta ,\eta_{\ast} ,g_{m})$ onto Riemannian
manifolds $(N^{\ast},g_{n})$. Then the ns conditions for the $(ker\Psi
_{\ast })$ to be completely geodesic foliation on $M^{\ast}$ are
\begin{eqnarray*}
g_{n}(\nabla _{\omega V}\Psi _{\ast }Z,\Psi _{\ast }\omega U)&
=\lambda_{\ast} ^{2}\{g_{m}(\beta T_{U}\phi V+A_{\omega V}\phi
U+g_{m}(\omega V,\omega U)grad\ ln\lambda_{\ast} \\
& -g_{m}(\phi U,V)\zeta ,Z)\},
\end{eqnarray*}%
and
\[
T_{U}\omega V+\hat{\nabla}_{U}\phi V\in \Gamma (D^{\ast}_{1}).
\]%
for $U,V\in \Gamma (ker\Psi _{\ast })$, $Z\in \Gamma (\mu )$ and $W\in
\Gamma (D^{\ast}_{2})$.
\end{th}

\noindent {\textbf Proof.} It is to be noted that the ns conditions for the $%
(ker\Psi _{\ast })$ to be completely geodesic foliation on $M^{\ast}$ are $%
g_{m}(\nabla _{U}^{m}V,Z)=0$ and $g_{m}(\nabla _{U}^{m}V,\varphi_{\ast} W)=0$
for $U,V\in \Gamma (ker\Psi _{\ast })$, $Z\in \Gamma (\mu )$ and $W\in
\Gamma (D^{\ast}_{2})$. Then with the help of (\ref{EQ2.1}), (\ref{EQ2.2}), (%
\ref{EQ2.4}) and (\ref{EQ3.3}), we infer
\begin{eqnarray*}
g_{m}(\nabla _{U}^{m}V,Z)& =g_{m}(\nabla _{U}^{m}\phi V,\varphi_{\ast}
Z)+g_{m}(\phi U,\nabla _{\omega V}^{m}Z)+g_{m}(\omega U,\nabla _{\omega
V}^{m}Z) \\
& +\eta_{\ast} (Z)g_{m}(\varphi_{\ast} U,V).
\end{eqnarray*}%
Since $\Psi $ is a sic $\zeta ^{\perp }$-Rs, from (\ref{EQ2.10}), (\ref%
{EQ2.13}), (\ref{EQ2.15}) and Lemma \ref{LE3}, we get
\begin{eqnarray*}
g_{m}(\nabla _{U}^{m}V,Z)& =g_{m}(T_{U}\phi V,\varphi_{\ast} Z)+g_{m}(\phi
U,A_{\omega V}Z) \\
&-\frac{1}{\lambda_{\ast} ^{2}}g_{m}(grad\ ln\lambda_{\ast} ,Z)g_{n}(\Psi
_{\ast }\omega V,\Psi _{\ast }\omega U) \\
& +\frac{1}{\lambda_{\ast} ^{2}}g_{n}(\nabla _{\omega V}^{\Psi }\Psi _{\ast
}Z,\Psi _{\ast }\omega U)+\eta_{\ast} (Z)g_{m}(\varphi_{\ast} U,V).
\end{eqnarray*}%
From the previous equation, we have
\begin{eqnarray}
g_{m}(\nabla _{U}^{m}V,Z)& =g_{m}(-\beta T_{U}\phi V-A_{\omega V}\phi
U-g_{m}(\omega V,\omega U)(grad\ ln\lambda_{\ast} ),Z)  \label{EQ3.17} \\
& +g_{m}(\phi U,V)\zeta ,Z)+\frac{1}{\lambda_{\ast} ^{2}}g_{n}(\nabla
_{\omega V}^{\Psi }\Psi _{\ast }Z,\Psi _{\ast }\omega U).  \nonumber
\end{eqnarray}%
Now, using (\ref{EQ2.1}), (\ref{EQ2.2}) (2.8), (2.9) and (\ref{EQ3.3}), we
infer
\begin{eqnarray}
g_{m}(\nabla _{U}^{m}V,\varphi_{\ast} W) &=&-g_{m}(\varphi_{\ast}
(T_{U}\omega V+\hat{\nabla}_{U}\phi V),\varphi_{\ast} W)  \nonumber \\
&=&-g_{m}(\omega (T_{U}\omega V+\hat{\nabla}_{U}\phi V),\varphi_{\ast} W)
\label{EQ3.18}
\end{eqnarray}%
So the results follows from (\ref{EQ3.17}) and (\ref{EQ3.18}).

Next, we offer certain conditions for dilation $\lambda_{\ast} $ is a
constant on $\mu $. For which, we need the following convenient definition.

\begin{defn}
Let $\Psi $ be a sic $\zeta ^{\perp }$-Rs from acm manifolds $%
(M,g_{m},J,\zeta ,\eta_{\ast} )$ onto Riemannian manifolds $(N^{\ast},g_{n})$%
. Then we call that $\mu $ is horizontal along $(ker\Psi _{\ast })$ if $%
\nabla _{U}^{m}Z\in \Gamma (\mu )$ for $Z\in \Gamma (\mu )$ and $U\in \Gamma
(ker\Psi _{\ast })$.
\end{defn}

\begin{cor}
Let $\Psi $ be a sic $\zeta ^{\perp }$-Rs from Sasakian manifolds $%
(M^{\ast},\varphi_{\ast} ,\zeta ,\eta_{\ast} ,g_{m})$ onto Riemannian
manifolds $(N^{\ast},g_{n})$ such that $\mu $ is horizontal along $(ker\Psi
_{\ast })$. Then the ns condition for the $\Psi $ to be constant on $\mu $
is
\begin{equation}
\lambda_{\ast} ^{2}g_{m}(\beta T_{U}\phi V+A_{\omega V}\phi U-g_{m}(\phi
U,V)\zeta ,Z)=g_{n}(\nabla _{\omega V}\Psi _{\ast }Z,\Psi _{\ast }\omega U),
\label{EQ3.19}
\end{equation}%
for $U,V\in (ker\Psi _{\ast })$ and $Z\in \Gamma (\mu )$.
\end{cor}

\noindent{\textbf Proof.} Let $U,V\in \Gamma (ker\Psi _{\ast })$, $Z\in
\Gamma (\mu )$. Then $g_{m}(\nabla _{U}^{m}V,Z)=-g_{m}(V,\nabla _{U}^{m}Z)$.
Since $\mu $ is horizontal along $(ker\Psi _{\ast })$, so $g_{m}(V,\nabla
_{U}^{m}Z)=0$.

Using (\ref{EQ3.17}) and (\ref{EQ3.19}), we infer
\begin{equation}
g_{m}(\omega V,\omega U)g_{m}(grad\ ln\lambda_{\ast} ,Z)=0.  \label{EQ3.20}
\end{equation}%
Replacing $U$ by $V$ in the foregoing equation, we may say that $%
\lambda_{\ast} $ is a constant on $\Gamma (\mu )$.\newline
The converse follows from (\ref{EQ3.17}).

\begin{th}
Let $\Psi $ be a sic $\zeta ^{\perp }$-Rs from Sasakian manifolds $%
(M^{\ast},\varphi_{\ast} ,\zeta ,\eta_{\ast} ,g_{m})$ onto Riemannian
manifolds $(N^{\ast},g_{n})$. Then the ns conditions for $M^{\ast}$ to be a
locally product manifold of the form $M_{\ker \Psi _{\ast }}\times
_{\lambda_{\ast} }M_{\ker \Psi _{\ast }^{\perp }}$ are
\begin{eqnarray*}
g_{n}(\nabla _{\omega V}\Psi _{\ast }Z,\Psi _{\ast }\omega U)&
=\lambda_{\ast} ^{2}\{g_{m}(\beta T_{U}\phi V+A_{\omega V}\phi
U+g_{m}(\omega V,\omega U)grad\ ln\lambda_{\ast} \\
& -g_{m}(\phi U,V)\zeta ,Z)\},
\end{eqnarray*}%
\[
T_{U}\omega V+\hat{\nabla}_{U}\phi V\in \Gamma (D^{\ast}_{1}).
\]%
and
\[
A_{P}\varphi_{\ast} Q+{\mathcal{V}_{1}}\nabla _{P}^{m}\alpha Q\in \Gamma
(D^{\ast}_{2}),
\]%
\begin{eqnarray*}
g_{n}(\nabla _{P}\Psi _{\ast }\varphi_{\ast} W,\Psi _{\ast }\beta Q)=&
\lambda_{\ast} ^{2}g_{m}(-A_{P}\alpha Q+\beta Q(ln\lambda_{\ast}
)P-g_{m}(P,\beta Q)({grad}\ln \lambda_{\ast} ) \\
& -P(ln\lambda_{\ast} )\beta Q+\eta_{\ast} (Q)P,\varphi_{\ast} W)
\end{eqnarray*}%
for $P,Q\in \Gamma (ker\Psi _{\ast })^{\perp }$, $U,V,W\in \Gamma (ker\Psi
_{\ast })$.
\end{th}

Let $\Psi $ be a sic $\zeta ^{\perp }$-Rs from Sasakian manifolds $%
(M^{\ast},\varphi_{\ast} ,\zeta ,\eta_{\ast} ,g_{m})$ onto Riemannian
manifolds $(N^{\ast},g_{n})$. For $U\in \Gamma (D^{\ast}_{2})$ and $Q\in
\Gamma (\mu )$, we have
\[
\frac{1}{\lambda_{\ast} ^{2}}g_{n}(\Psi _{\ast }\varphi_{\ast} U,\Psi _{\ast
}Q)=g_{m}(\varphi_{\ast} U,Q)=0,
\]%
which shows that the distributions $\Psi _{\ast }(\varphi_{\ast}
D^{\ast}_{2})$ and $\Psi _{\ast }(\mu )$ are orthogonal.

Now, we explore the geometric properties of the leaves of the distributions $%
D^{\ast}_{1}$ and $D^{\ast}_{2}$.

\begin{th}
$\label{th6.}$Let $\Psi $ be a sic $\zeta ^{\perp }$-Rs from Sasakian
manifolds $(M^{\ast},\varphi_{\ast} ,\zeta ,\eta_{\ast} ,g_{m}) $ onto
Riemannian manifolds $(N^{\ast},g_{n})$. Then the ns conditions for the
distribution $D^{\ast}_{1}$ to be completely geodesic foliations on $M$ are
\[
(\nabla \Psi _{\ast })(U,\varphi_{\ast} V)\in \Gamma (\Psi _{\ast }\mu )
\]%
and
\[
\frac{1}{\lambda_{\ast} ^{2}}g_{n}(\nabla \Psi _{\ast })(U,\varphi_{\ast}
V),\Psi _{\ast }\beta Z)=g_{m}(V,T_{U}\omega \alpha Z+\hat{\nabla}_{U}\phi
\alpha Z)-\eta_{\ast} (\beta Z)g_{m}(U,V)
\]%
for $U,V\in \Gamma (D^{\ast}_{1})$ and $Z\in \Gamma (ker\Psi _{\ast
})^{\perp }$.
\end{th}

\noindent {\textbf Proof.} The ns conditions for the distribution $%
D^{\ast}_{1}$ to be completely geodesic foliations on $M$ are $g_{m}(\nabla
_{U}^{m}V,W)=0$ and $g_{m}(\nabla _{U}^{m}V,Z)=0$ for $U,V\in \Gamma
(D^{\ast}_{1})$, $W\in \Gamma (D^{\ast}_{2})$ and $Z\in \Gamma (\ker \Psi
_{\ast })^{\perp }$. Using (\ref{EQ2.1}), (\ref{EQ2.2}), (\ref{EQ2.4}) and (%
\ref{EQ3.3}), we obtain
\begin{eqnarray*}
g_{m}(\nabla _{U}^{m}V,W)& =g_{m}(\varphi_{\ast} \nabla
_{U}^{m}V,\varphi_{\ast} W)-\eta_{\ast} (\nabla _{U}^{m}V)\eta_{\ast} (W) \\
& =g_{m}(\nabla _{U}^{m}\varphi_{\ast} V,\varphi_{\ast} W)
\end{eqnarray*}%
From (\ref{EQ2.15}) and using the fact that $\Psi $ is a sic $\zeta ^{\perp }
$-Rs, we have
\begin{equation}
g_{m}(\nabla _{U}^{m}V,W))=-\frac{1}{\lambda_{\ast} ^{2}}g_{n}((\nabla \Psi
_{\ast })(U,\varphi_{\ast} V),\Psi _{\ast }\varphi_{\ast} W)  \label{EQ3.21}
\end{equation}%
Using (\ref{EQ2.1}), (\ref{EQ2.2}), (\ref{EQ2.4}), (\ref{EQ2.10}) and (\ref%
{EQ3.4}), we have
\begin{eqnarray*}
g_{m}(\nabla _{U}^{m}V,Z) &=&-g_{m}(V,\nabla _{U}^{m}Z)=g_{m}(V,\nabla
_{U}^{m}\varphi_{\ast} ^{2}Z) \\
&=&g_{m}(V,\nabla _{U}^{m}\varphi_{\ast} \alpha Z)+g_{m}(V,\nabla
_{U}^{m}\varphi_{\ast} \beta Z) \\
&=&g_{m}(V,\nabla _{U}^{m}\varphi_{\ast} \alpha Z)+g_{m}(\nabla
_{U}^{m}\varphi_{\ast} V,\beta Z)-\eta_{\ast} (\beta Z)g_{m}(U,V)
\end{eqnarray*}%
Now, using (\ref{EQ2.11}), (\ref{EQ2.15}) and (\ref{EQ3.3}), we get
\begin{eqnarray}
g_{m}(\nabla _{U}^{m}V,Z) &=g_{m}(V,T_{U}\omega \alpha Z)+g_{m}(V,\hat{\nabla%
}_{U}\phi \alpha Z)  \nonumber \\
&-\frac{1}{\lambda_{\ast} ^{2}}g_{m}((\nabla \Psi _{\ast })(U,\varphi_{\ast}
V),\Psi _{\ast }\beta Z)  \nonumber \\
&-\eta_{\ast} (\beta Z)g_{m}(U,V).  \label{EQ3.22}
\end{eqnarray}%
The proof of the theorem comes from (\ref{EQ3.21}) and (\ref{EQ3.22}).

\begin{th}
\label{th7}Let $\Psi $ be a sic $\zeta ^{\perp }$-Rs from Sasakian manifolds
$(M^{\ast},\varphi_{\ast} ,\zeta ,\eta_{\ast} ,g_{m}) $ onto Riemannian
manifolds $(N^{\ast},g_{n})$. Then the ns conditions for the distribution $%
D^{\ast}_{2}$ to be a completely geodesic foliations on $M^{\ast}$ are
\[
(\nabla \Psi _{\ast })(W,\varphi_{\ast} U)\in \Gamma (\Psi _{\ast }\mu )
\]%
and
\begin{eqnarray*}
\frac{1}{\lambda_{\ast} ^{2}}g_{n}(\nabla _{\varphi_{\ast} \xi}^{\Psi }\Psi
_{\ast }\varphi_{\ast} W,\Psi _{\ast }\varphi_{\ast} \beta Z)
&=&g_{m}(\varphi_{\ast} \xi,T_{W}\alpha Z)-g_{m}(W,\xi)g_{m}({\mathcal{H}_{1}%
}grad\ ln\lambda_{\ast} ,\varphi_{\ast} \beta Z) \\
&&+\eta_{\ast} (\beta Z)g_{m}(W,\xi).
\end{eqnarray*}%
for $W,\xi\in \Gamma (D^{\ast}_{2})$, $U\in \Gamma (D^{\ast}_{1})$ and $Z\in
\Gamma (ker\Psi _{\ast }^{\perp })$.
\end{th}

\noindent {\textbf Proof.} The ns conditions for the distribution $%
D^{\ast}_{2}$ to be a completely geodesic foliations on $M^{\ast}$ are $%
g_{m}(\nabla _{W}^{m}\xi,U)=0$ and $g_{m}(\nabla _{W}^{m}\xi,Z)=0$ for $%
W,\xi\in \Gamma (D^{\ast}_{2})$, $U\in \Gamma (D^{\ast}_{1})$ and $Z\in
\Gamma ((\ker \Psi _{\ast })^{\perp })$. Using (\ref{EQ2.1}), (\ref{EQ2.2}),
(\ref{EQ2.4}) and (\ref{EQ3.3}), we come to the next equation
\begin{eqnarray*}
g_{m}(\nabla _{W}^{m}\xi,U))& =g_{m}(\varphi_{\ast} \nabla
_{W}^{m}\xi,\varphi_{\ast} U)+\eta_{\ast} (\nabla _{W}^{m}\xi)\eta_{\ast} (U)
\\
& =g_{m}(\nabla _{W}^{m}\varphi_{\ast} \xi,\varphi_{\ast}
U)=-g_{m}(\varphi_{\ast} \xi,\nabla _{W}^{m}\varphi_{\ast} U)
\end{eqnarray*}%
From (\ref{EQ2.15}) and definition of sic $\zeta ^{\perp } $-Rs, we get
\begin{equation}
g_{m}(\nabla _{W}^{m}\xi,U))=\frac{1}{\lambda_{\ast} ^{2}}g_{m}((\nabla \Psi
_{\ast })(W,\varphi_{\ast} U),\Psi _{\ast }\varphi_{\ast} \xi).
\label{EQ3.23}
\end{equation}%
With the help of (\ref{EQ2.1}), (\ref{EQ2.2}), (\ref{EQ2.4}), (\ref{EQ2.10})
and (\ref{EQ3.4}), we obtain
\[
g_{m}(\nabla _{W}^{m}\xi,Z))=-g_{m}(\varphi_{\ast} \xi,T_{W}\alpha
Z)+g_{m}(\nabla _{\varphi_{\ast} \xi}^{m}W,\beta Z).
\]%
\[
g_{m}(\nabla _{W}^{m}\xi,Z)=-g_{m}(\varphi_{\ast} \xi,T_{W}\alpha
Z)+g_{m}(\nabla _{\varphi_{\ast} \xi}^{m}\varphi_{\ast} W,\varphi_{\ast}
\beta Z)-\eta_{\ast} (\beta Z)g_{m}(W,\xi).
\]%
Then from (\ref{EQ2.15}), (\ref{EQ3.4}) and Lemma \ref{LE3}, we obtain
\begin{eqnarray}
g_{m}(\nabla _{W}^{m}\xi,Z)& =-g_{m}(\varphi_{\ast} \xi,T_{W}\alpha
Z)+g_{m}(W,\xi)g_{m}({\mathcal{H}_{1}}grad\ ln\lambda_{\ast} ,\varphi_{\ast}
\beta Z)  \label{EQ3.24} \\
& +\frac{1}{\lambda_{\ast} ^{2}}g_{n}(\nabla _{\varphi_{\ast} \xi}^{\Psi
}\Psi _{\ast }\varphi_{\ast} W,\Psi _{\ast }\varphi_{\ast} \beta
Z)-\eta_{\ast} (\beta Z)g_{m}(W,\xi).  \nonumber
\end{eqnarray}%
The proof of the theorem comes from (\ref{EQ3.23}) and (\ref{EQ3.24}).

From Theorem \ref{th6.} and \ref{th7}, we have the following theorem:

\begin{th}
Let $\Psi $ be a sic $\zeta ^{\perp }$-Rs from Sasakian manifolds $%
(M^{\ast},\varphi_{\ast} ,\zeta ,\eta_{\ast} ,g_{m})$ onto Riemannian
manifolds $(N^{\ast},g_{n})$. Then the ns conditions for the fibers of $\Psi
$ to be locally product manifold are
\[
(\nabla \Psi _{\ast })(U,\varphi_{\ast} V)\in \Gamma (\Psi _{\ast }\mu )
\]%
and
\[
(\nabla \Psi _{\ast })(W,\varphi_{\ast} U)\in \Gamma (\Psi _{\ast }\mu )
\]%
for any $U,V\in \Gamma (D^{\ast}_{1})$ and $W\in \Gamma (D^{\ast}_{2})$.
\end{th}

\section{Harmonicity of sic $\protect\zeta ^{\perp }$-RS \label{sect-har}}

In this part, we provide the ns conditions for a sic $\zeta ^{\perp }$-Rs to
be harmonic. We also look into the ns conditions for such submersions to be
completely geodesic. By decomposition of all over space of sic $\zeta
^{\perp }$-Rs, we have the following Lemma.

\begin{lem}
\label{LE4.1} Let $\Psi :(M^{2(p+q+r)+1},g_{m},\varphi_{\ast} )\rightarrow
(N^{q+2r+1},g_{n})$ be a sic $\zeta ^{\perp }$-Rs from acm manifolds $%
(M^{\ast},\varphi_{\ast} ,\zeta ,\eta_{\ast} ,g_{m})$ onto Riemannian
manifolds $(N^{\ast},g_{n})$. Then the tension field $\tau $ of $\Psi $ is
\begin{equation}
\tau (\Psi )=-(2p+q)\Psi _{{\ast }}(\mu ^{\ker \Psi _{{\ast }%
}})+(2-q-2r)\Psi _{{\ast }}({grad}\ln \lambda_{\ast} ),  \label{EQ4.1}
\end{equation}%
where $\mu ^{\ker \Psi _{{\ast }}}$ is the mean curvature vector field of
the distribution of $\ker \Psi _{{\ast }}$.
\end{lem}

\noindent {\textbf Proof:} It is easy to prove this lemma with the help of
the method \cite{Akyol}.

In view of Lemma \ref{LE4.1}, we have the following result.

\begin{th}
Let $\Psi :(M^{2(p+q+r)+1},g_{m},\varphi_{\ast} )\rightarrow
(N^{q+2r+1},g_{n})$ be a sic $\zeta ^{\perp }$-Rs from acm manifolds $%
(M^{\ast},\varphi_{\ast} ,\zeta ,\eta_{\ast} ,g_{m})$ onto Riemannian
manifolds $(N^{\ast},g_{n})$ such that $q+2r\neq 2$. Then any two of the
conditions imply the third:\newline

\indent(i) $\Psi$ is harmonic\newline

\indent(ii) The fibers are minimal.\newline

\indent(iii) $\Psi$ is a horizontally homothetic map.
\end{th}

\noindent {\textbf Proof:} Using (\ref{EQ4.1}), we can easily find the
results.

\begin{cor}
Let $\Psi :(M^{2(p+q+r)+1},g_{m},\varphi_{\ast} )\rightarrow
(N^{q+2r+1},g_{n})$ be a sic $\zeta ^{\perp }$-Rs from acm manifolds $%
(M^{\ast},\varphi_{\ast} ,\zeta ,\eta_{\ast} ,g_{m})$ onto Riemannian
manifolds $(N^{\ast},g_{n})$. If $q+2r=2$ then the ns condition for the $%
\Psi $ to be harmonic is the minimal fiber.
\end{cor}

A map is claimed to be totally geodesic map if it maps every geodesic in the
total manifold into a geodesic within the base manifold in proportion to arc
lengths. We recall that a differentiable map $\Psi $ between two Riemannian
manifolds is called completely geodesic if $(\nabla \Psi _{\ast })(U,V)=0\ \
\forall U,V\in \Gamma (TM^{\ast})$.

We now, present the subsequent definition:

\begin{defn}
Let $\Psi $ be a sic $\zeta ^{\perp }$-Rs from acm manifolds $%
(M^{\ast},\varphi_{\ast} ,\zeta ,\eta_{\ast} ,g_{m})$ onto Riemannian
manifolds $(N^{\ast},g_{n})$. Then $\Psi $ is called a $(\varphi_{\ast}
D^{\ast}_{2},\mu )$-completely geodesic map if
\[
(\nabla \Psi _{\ast })(\varphi_{\ast} W,Q)=0\;\ {\rm for}\ W\in \Gamma
(D^{\ast}_{2})\ {\rm and}\ Q\in \Gamma (\mu ).
\]
\end{defn}

In this sequence, we have the following important result.

\begin{th}
Let $\Psi$ be a sic $\zeta ^{\perp }$-Rs from acm manifolds $%
(M^{\ast},\varphi_{\ast} ,\zeta ,\eta_{\ast} ,g_{m})$ onto Riemannian
manifolds $(N^{\ast},g_{n})$. The ns condition for $\Psi$ to be a $%
(\varphi_{\ast} D^{\ast}_{2},\mu )$-completely geodesic map if $\Psi $ is a
horizontally homothetic map.
\end{th}

\noindent {\textbf Proof.} For $W\in \Gamma (D^{\ast}_{2})$ and $Q\in \Gamma
(\mu )$, from Lemma \ref{LE3}, we infer
\[
(\nabla \Psi _{\ast })(\varphi_{\ast} W,Q)=\varphi_{\ast} W(\ln
\lambda_{\ast} )\Psi _{\ast }Q+Q(ln\lambda_{\ast} )\Psi _{\ast
}\varphi_{\ast} W-g_{m}(\varphi_{\ast} W,Q)\Psi _{\ast }({grad}\ln
\lambda_{\ast} ).
\]%
Since $\Psi $ is a horizontally homothetic therefore $(\nabla \Psi _{\ast
})(\varphi_{\ast} W,Q)=0$.\newline
Conversely, if $(\nabla \Psi _{\ast })(\varphi_{\ast} W,Q)=0,$ we have
\begin{equation}
\varphi_{\ast} W(\ln \lambda_{\ast} )\Psi _{\ast }Q+Q(\ln \lambda_{\ast}
)\Psi _{\ast }\varphi_{\ast} W=0  \label{EQ4.2}
\end{equation}%
Taking scaler product in (\ref{EQ4.2}) with $\Psi _{\ast }\varphi_{\ast} W$
and since $\Psi $ is a sic $\zeta ^{\perp }$-Rs, we infer
\[
g_{m}({grad}\ln \lambda_{\ast} ,\varphi_{\ast} W)g_{n}(\Psi _{\ast }Q,\Psi
_{\ast }\varphi_{\ast} W)+g_{m}({grad}\ln \lambda_{\ast} ,Q)g_{n}(\Psi
_{\ast }\varphi_{\ast} W,\Psi _{\ast }\varphi_{\ast} W)=0.
\]%
Since $||\varphi_{\ast} W||\neq 0$, then the above equation reduces to
\[
g_{m}({\mathcal{H}_{1}}{grad}\ln \lambda_{\ast} ,Q)=0,
\]%
which implies that $\lambda_{\ast} $ is a constant on $\Gamma (\mu )$.
\newline
Taking scaler product in (\ref{EQ4.2}) with $\Psi _{\ast }Q$ ,we get
\[
g_{m}({grad}\ln \lambda_{\ast} ,\varphi_{\ast} W)g_{n}(\Psi _{\ast }Q,\Psi
_{\ast }Q)+g_{m}({grad}\ln \lambda_{\ast} ,Q)g_{n}(\Psi _{\ast
}\varphi_{\ast} W,\Psi _{\ast }Q)=0.
\]%
Since $||Q||\neq 0$, then the above equation reduces to
\[
g_{m}({grad}\ln \lambda_{\ast} ,\varphi_{\ast} W)=0,
\]%
it follows that $\lambda_{\ast} $ is a constant on $\Gamma
(\varphi_{\ast} D^{\ast}_{2})$. Thus $\lambda_{\ast} $ is a constant
on $\Gamma ((ker\Psi _{\ast })^{\perp })$, which over the
proof.\newline At last, we infer the ns conditions for sic $\zeta
^{\perp }$-Rs to be completely geodesic.

\begin{th}
Let $\Psi $ be a sic $\zeta ^{\perp }$-Rs from Sasakian manifolds $%
(M^{\ast},\varphi_{\ast} ,\zeta ,\eta_{\ast} ,g_{m})$ onto Riemannian
manifolds $(N^{\ast},g_{n})$. Then the ns conditions for the $\Psi $ to be
completely geodesic map are\newline

\indent(a) $\beta T_{U}\varphi_{\ast} V+\omega \hat{\nabla}%
_{U}\varphi_{\ast} V-\eta_{\ast} (T_{U}V)\zeta =0,\ U,V\in \Gamma
(D^{\ast}_{1})$,\newline

\indent(b) $\omega T_{U}\varphi_{\ast} W+\beta H\hat{\nabla}%
_{U}\varphi_{\ast} W-\eta_{\ast} (T_{U}^{m}W)\zeta =0,\ U\in \Gamma (ker\Psi
_{\ast }),\ W\in \Gamma (D^{\ast}_{2})$,\newline

\indent(c) $\Psi $ is a horizontally homothetic map,\newline

\indent (d) $\omega {\mathcal{V}_{1}}\nabla _{P}\phi U+\beta A_{P}\phi
U+\beta {\mathcal{H}_{1}}\nabla _{P}\omega U+\omega A_{P}\omega
U-\eta_{\ast} \left( A_{P}U\right) \zeta =0$, $U\in \Gamma (ker\Psi _{\ast })
$, $P\in \Gamma ((ker\Psi _{\ast })^{\perp })$.
\end{th}

\noindent {\textbf Proof.} (a) For $U,V\in \Gamma (D_{1}^{\ast })$, using (%
\ref{EQ2.1}), (\ref{EQ2.2}), (\ref{EQ2.4}), (\ref{EQ2.10}) and (\ref{EQ2.15}%
), we have
\begin{eqnarray*}
(\nabla {\Psi _{\ast }})(U,V)& =\nabla _{U}\Psi _{\ast }V-\Psi _{\ast
}(\nabla _{U}^{m}V) \\
& =\Psi _{\ast }(\varphi _{\ast }^{2}\nabla _{U}^{m}V-\eta _{\ast }(\nabla
_{U}^{m}V)\zeta ) \\
& =\Psi _{\ast }(\varphi _{\ast }\nabla _{U}^{m}\varphi _{\ast }V-\eta
_{\ast }(\nabla _{U}^{m}V)\zeta ) \\
& =\Psi _{\ast }(\varphi _{\ast }(T_{U}\varphi _{\ast }V+\hat{\nabla}%
_{U}\varphi _{\ast }V)-\eta _{\ast }(T_{U}V)\zeta ).
\end{eqnarray*}%
With the help of (\ref{EQ3.3}) and (\ref{EQ3.4}) in the above expression
turns into
\[
(\nabla {\Psi _{\ast }})(U,V)=\Psi _{\ast }(\alpha T_{U}\varphi _{\ast
}V+\beta T_{U}\varphi _{\ast }V+\phi \hat{\nabla}_{U}\varphi _{\ast
}V+\omega \hat{\nabla}_{U}\varphi _{\ast }V-\eta _{\ast }(\nabla
_{U}^{m}V)\zeta ),
\]%
since $\alpha T_{U}\varphi _{\ast }V+\phi \hat{\nabla}_{U}\varphi _{\ast
}V\in \Gamma ((ker\Psi _{\ast })$, we infer
\[
(\nabla {\Psi _{\ast }})(U,V)=\Psi _{\ast }(\beta T_{U}\varphi _{\ast
}V+\omega \hat{\nabla}_{U}\varphi _{\ast }V-\eta _{\ast }(\nabla
_{U}^{m}V)\zeta ).
\]%
Since $\Psi $ is a linear isomorphism between $(ker\Psi _{\ast })^{\perp }$
and $\Gamma (TN)$, so $(\nabla {\Psi _{\ast }})(U,V)=0$ if and only if $%
\beta T_{U}\varphi _{\ast }V+\omega \hat{\nabla}_{U}\varphi _{\ast }V-\eta
_{\ast }(\nabla _{U}^{m}V)\zeta =0.$\newline
(b) For $U\in \Gamma (ker\Psi _{\ast })$ and $W\in \Gamma (D_{2}^{\ast })$,
using (\ref{EQ2.1}), (\ref{EQ2.2}), (\ref{EQ2.4}) and \ref{EQ2.15}), we come
to the next equation
\begin{eqnarray*}
(\nabla {\Psi _{\ast }})(U,W)& =\nabla _{U}\Psi _{\ast }W-\Psi _{\ast
}(\nabla _{U}^{m}W) \\
& =\Psi _{\ast }(\varphi _{\ast }^{2}\nabla _{U}^{m}W-\eta _{\ast }(\nabla
_{U}^{m}W)\zeta ) \\
& =\Psi _{\ast }(\varphi _{\ast }\nabla _{U}^{m}\varphi _{\ast }W-\eta
_{\ast }(T_{U}^{m}W)\zeta ) \\
& =\Psi _{\ast }(\varphi _{\ast }(T_{U}\varphi _{\ast }W+{\mathcal{H}_{1}}%
\nabla _{U}\varphi _{\ast }W)-\eta _{\ast }(T_{U}^{m}W)\zeta ).
\end{eqnarray*}%
In view of (\ref{EQ3.3}) and (\ref{EQ3.4}) in the foregoing equation, we
gain
\[
(\nabla {\Psi _{\ast }})(U,W)=\Psi _{\ast }(\phi T_{U}\varphi _{\ast
}W+\omega T_{U}\varphi _{\ast }W+\alpha {\mathcal{H}_{1}}\nabla _{U}\varphi
_{\ast }W+\beta {\mathcal{H}_{1}}\nabla _{U}\varphi _{\ast }W-\eta _{\ast
}(T_{U}^{m}W)\zeta )
\]%
Since $\phi T_{U}\varphi _{\ast }W+\alpha {\mathcal{H}_{1}}\nabla
_{U}\varphi _{\ast }W\in \Gamma (\ker \Psi _{\ast })$, we infer
\[
(\nabla {\Psi _{\ast }})(U,W)=\Psi _{\ast }(\omega T_{U}\varphi _{\ast
}W+\beta {\mathcal{H}_{1}}\nabla _{U}\varphi _{\ast }W-\eta _{\ast
}(T_{U}^{m}W)\zeta ).
\]%
Since $\Psi $ is a linear isomorphism between $(ker\Psi _{\ast })^{\perp }$
and $\Gamma (TN)$, so $(\nabla {\Psi _{\ast }})(U,W)=0$ if and only if $%
\omega T_{U}\varphi _{\ast }W+\beta H\hat{\nabla}_{U}\varphi _{\ast }W-\eta
_{\ast }(T_{U}^{m}W)\zeta =0.$\newline
(c) For $U,V\in \Gamma (\mu )$, from Lemma \ref{LE3}, we get
\[
(\nabla \Psi _{\ast })(U,V)=U(\ln \lambda _{\ast })\Psi _{\ast }V+V(\ln
\lambda _{\ast })\Psi _{\ast }U-g_{m}(U,V)\Psi _{\ast }({grad}\ln \lambda
_{\ast }).
\]%
In the above equation, replacing $V$ by $\varphi _{\ast }U$, we infer
\begin{eqnarray*}
(\nabla \Psi _{\ast })(U,\varphi _{\ast }U)& =U(\ln \lambda _{\ast })\Psi
_{\ast }\varphi _{\ast }U+\varphi _{\ast }U(\ln \lambda _{\ast })\Psi _{\ast
}U-g_{m}(U,\varphi _{\ast }U)\Psi _{\ast }({grad}\ln \lambda _{\ast }). \\
& =U(\ln \lambda _{\ast })\Psi _{\ast }\varphi _{\ast }U+\varphi _{\ast
}U(\ln \lambda _{\ast })\Psi _{\ast }U.
\end{eqnarray*}%
If $(\nabla \Psi _{\ast })(U,\varphi _{\ast }U)=0$, we have
\begin{equation}
U(\ln \lambda _{\ast })\Psi _{\ast }\varphi _{\ast }U+\varphi _{\ast }U(\ln
\lambda _{\ast })\Psi _{\ast }U=0.  \label{EQ4.3}
\end{equation}%
Taking scaler product of the above equation with $\Psi _{\ast }\varphi
_{\ast }U$, we infer
\[
g_{m}({grad}\ln \lambda _{\ast },U)g_{n}(\Psi _{\ast }\varphi _{\ast }U,\Psi
_{\ast }\varphi _{\ast }U)+g_{m}(gradln\lambda _{\ast },\varphi _{\ast
}U)g_{n}(\Psi _{\ast }U,\Psi _{\ast }\varphi _{\ast }U)=0,
\]%
which shows that $\lambda _{\ast }$ is a constant $\Gamma (\mu )$.\newline
In a similar pattern, for $U,V\in \Gamma (D_{2}^{\ast })$, using Lemma \ref%
{LE3}, we get
\begin{eqnarray*}
(\nabla \Psi _{\ast })(\varphi _{\ast }U,\varphi _{\ast }V)& =\varphi _{\ast
}U(ln\lambda _{\ast })\Psi _{\ast }\varphi _{\ast }V+\varphi _{\ast
}V(ln\lambda _{\ast })\Psi _{\ast }\varphi _{\ast }U \\
& -g_{m}(\varphi _{\ast }U,\varphi _{\ast }V)\Psi _{\ast }(grad\ ln\lambda
_{\ast }).
\end{eqnarray*}%
In the above equation, replacing $V$ by $U$, we get
\begin{eqnarray}
(\nabla \Psi _{\ast })(\varphi _{\ast }U,\varphi _{\ast }U)& =\varphi _{\ast
}U(ln\lambda _{\ast })\Psi _{\ast }\varphi _{\ast }U+\varphi _{\ast
}U(ln\lambda _{\ast })\Psi _{\ast }\varphi _{\ast }U  \nonumber \\
& -g_{m}(\varphi _{\ast }U,\varphi _{\ast }U)\Psi _{\ast }(grad\ ln\lambda
_{\ast })  \nonumber  \label{EQ4.4} \\
& =2\varphi _{\ast }U(ln\lambda _{\ast })\Psi _{\ast }\varphi _{\ast
}U-g_{m}(\varphi _{\ast }U,\varphi _{\ast }U)\Psi _{\ast }(grad\ ln\lambda
_{\ast }).
\end{eqnarray}%
Consider $(\nabla \Psi _{\ast })(\varphi _{\ast }U,\varphi _{\ast }U)=0$ and
taking scalar product of above equation with $\Psi _{\ast }\varphi _{\ast }U$
and using the fact that $\Psi $ is sic $\zeta ^{\perp }$-Rs, we infer
\[
2g_{m}(grad\ ln\lambda _{\ast },\varphi _{\ast }U)g_{n}(\Psi _{\ast }\varphi
_{\ast }U,\Psi _{\ast }\varphi _{\ast }U)-g_{m}(\varphi _{\ast }U,\varphi
_{\ast }U)g_{n}(\Psi _{\ast }(grad\ ln\lambda _{\ast }),\Psi _{\ast }\varphi
_{\ast }U).
\]%
From the above expression, it follows that $\lambda _{\ast }$ is a constant
on $\Gamma (\varphi _{\ast }D_{2}^{\ast })$.\newline
Thus $\lambda _{\ast }$ is a constant on $\Gamma ((ker\Psi _{\ast })^{\perp
})$.

If $\Psi $ is a horizontally homothetic map, then ${\mathcal{H}_{1}}{grad}\
ln\lambda_{\ast} $ vanishes.\newline
Therefore $(\nabla \Psi _{\ast })(U,V)=0$ for $U,V\in \Gamma ((ker\Psi
_{\ast })^{\perp })$.\newline
(d) For $U\in \Gamma (ker\Psi _{\ast })$ and $P\in \Gamma ((ker\Psi _{\ast
})^{\perp })$, using (\ref{EQ2.1}), (\ref{EQ2.2}), (\ref{EQ2.4}), (\ref%
{EQ2.11}) and (\ref{EQ2.15}), we acquire
\begin{eqnarray*}
(\nabla {\Psi _{\ast }})(P,U)& =\nabla _{P}\Psi _{\ast }U-\Psi _{\ast
}(\nabla _{P}^{m}U) \\
& =\Psi _{\ast }(\varphi_{\ast} ^{2}\nabla _{P}^{m}U-\eta_{\ast} (\nabla
_{P}^{m}U)\zeta ) \\
& =\Psi _{\ast }(\varphi_{\ast} \nabla _{P}^{m}\varphi_{\ast} U-\eta_{\ast}
(\nabla _{P}^{m}U)\zeta ).
\end{eqnarray*}%
Then from (\ref{EQ2.10}), (\ref{EQ2.11}) and (\ref{EQ3.4}), we arrive at
\begin{eqnarray*}
(\nabla {\Psi _{\ast }})(P,U) &=&\Psi _{\ast }(\varphi_{\ast} \nabla
_{P}^{m}(\phi +\omega )U-\eta_{\ast} \left( A_{P}U+{\mathcal{V}_{1}}\nabla
_{P}U\right) \zeta ) \\
&=&\Psi _{\ast }(\varphi_{\ast} (A_{P}\phi U+{\mathcal{V}_{1}}\nabla
_{P}\phi U)+\varphi_{\ast} ({\mathcal{H}_{1}}\nabla _{P}\omega U+A_{P}\omega
U)) \\
&&-\Psi _{\ast }(\eta_{\ast} \left( A_{P}U\right) \zeta ).
\end{eqnarray*}%
Using (\ref{EQ3.3})and (\ref{EQ3.4}) in the above equation, we have
\begin{eqnarray*}
(\nabla {\Psi _{\ast }})(P,U) &=&\Psi _{\ast }(\phi {\mathcal{V}_{1}}\nabla
_{P}\phi U+\omega {\mathcal{V}_{1}}\nabla _{P}\phi U+\alpha A_{P}\phi
U+\beta A_{P}\phi U \\
&&+\alpha {\mathcal{H}_{1}}\nabla _{P}\omega U+\beta {\mathcal{H}_{1}}\nabla
_{P}\omega U+\phi A_{P}\omega U+\omega A_{P}\omega U-\eta_{\ast} \left(
A_{P}U\right) \zeta )
\end{eqnarray*}%
Since $\phi {\mathcal{V}_{1}}\nabla _{P}\phi U+\alpha A_{P}\phi U+\alpha {%
\mathcal{H}_{1}}\nabla _{P}\omega U+\phi A_{P}\omega U\in \Gamma (ker\Psi
_{\ast })$, we infer%
\[
(\nabla {\Psi _{\ast }})(P,U)=\Psi _{\ast }(\omega {\mathcal{V}_{1}}\nabla
_{P}\phi U+\beta A_{P}\phi U+\beta {\mathcal{H}_{1}}\nabla _{P}\omega
U+\omega A_{P}\omega U-\eta_{\ast} \left( A_{P}U\right) \zeta ).
\]%
Since $\Psi $ is a linear isomorphism between $(ker\Psi _{\ast })^{\perp }$
and $\Gamma(TN)$, so $(\nabla {\Psi _{\ast }})(P,U)=0$ if and only if $%
\omega {\mathcal{V}_{1}}\nabla _{P}\phi U+\beta A_{P}\phi U+\beta {\mathcal{H%
}_{1}}\nabla _{P}\omega U+\omega A_{P}\omega U-\eta_{\ast} \left(
A_{P}U\right) \zeta =0$.

\medskip

\noindent {\textbf Conflict of interest:} The authors declare that they have
no conflict of interest.

\noindent Department of Mathematics \\ University of
Calcutta
\\ 35 B.C. Road 700019, Kolkata\\ West Bengal, India.\\
Email: uc$\_$de@yahoo.com

\medskip
\noindent Department of Mathematics and Astronomy \\ University of Lucknow
\\ Lucknow, Uttar Pradesh, 226007, India.\\
Email: shashi.royal.lko@gmail.com

\medskip
\noindent Department of Mathematics and Statistics\\Dr.
Harisingh Gour University \\ Sagar-470003, Madhya Pradesh, India.\\
Email: pgupta@dhsgsu.edu.in
\end{document}